\documentclass[11pt]{article}
\usepackage{amssymb,amsmath,bm}
\usepackage{mathrsfs}
\usepackage{hyperref}
\usepackage{amsthm}
\hypersetup{
    colorlinks,
    citecolor=black,
    filecolor=black,
    linkcolor=black,
    urlcolor=black
}
\usepackage{constants}
\usepackage{enumerate}
\usepackage{tikz}
\begin{document}
\newcommand{\beas}{\begin{eqnarray*}}
\newcommand{\enas}{\end{eqnarray*}}
\newcommand{\bea}{\begin{eqnarray}}
\newcommand{\ena}{\end{eqnarray}}
\newcommand{\bms}{\begin{multline*}}
\newcommand{\ems}{\end{multline*}}
\newcommand{\bels}{\begin{align*}}
\newcommand{\enls}{\end{align*}}
\newcommand{\bel}{\begin{align}}
\newcommand{\enl}{\end{align}}
\newcommand{\qmq}[1]{\quad \mbox{#1} \quad}
\newcommand{\qm}[1]{\quad \mbox{#1}}
\newcommand{\nn}{\nonumber}
\newcommand{\bbox}{\hfill $\Box$}
\newcommand{\ignore}[1]{}
\newcommand{\tr}{\mbox{tr}}
\newcommand{\Bvert}{\left\vert\vphantom{\frac{1}{1}}\right.}
\newtheorem{theorem}{Theorem}[section]
\newtheorem{corollary}{Corollary}[section]
\newtheorem{conjecture}{Conjecture}[section]
\newtheorem{proposition}{Proposition}[section]
\newtheorem{remark}{Remark}[section]
\newtheorem{lemma}{Lemma}[section]
\newtheorem{definition}{Definition}[section]
\newtheorem{example}{Example}[section]
\newtheorem{condition}{Condition}[section]

\newcommand{\pf}{\noindent {\bf Proof:} }
\def\blfootnote{\xdef\@thefnmark{}\@footnotetext}
\newcommand{\lcolor}[1]{\textcolor{magenta}{#1}}
\newcommand{\ycolor}[1]{\textcolor{blue}{#1}}

\newcommand{\lcomm}[1]{\marginpar{\tiny\lcolor{#1}}}
\newcommand{\ycomm}[1]{\marginpar{\tiny\ycolor{#1}}}

\newcommand{\spa}{S}
\newcommand{\bbs}{\mathbb{S}}

\newcommand{\expect}[1]{E{\left[#1\right]}}

\newcommand{\tbbox}{\small \framebox[1.45\width]{$t$}}
\newcommand{\stbbox}{\small \framebox[1.45\width]{$st$}}
\newcommand{\sobbox}{\small \framebox[1.45\width]{$so$}}
\newcommand{\oobbox}{\small \framebox[1.45\width]{$oo$}}
\newcommand{\eebbox}{\small \framebox[1.45\width]{$11$}}

\newcommand{\qex}{\hfill \ensuremath{\blacktriangle}}
\newcommand{\xbbox}{\framebox[1.45\width]{$\bf\approx$}}

\title{A BKR operation for events occurring for disjoint reasons with high probability}
\author{Larry Goldstein and Yosef Rinott\\University of Southern California, and the Hebrew University of Jerusalem}

\maketitle

  \begin{abstract}
	Given events $A$ and $B$ on a product space $S=\prod_{i=1}^n S_i$, the set 
	$A \Box B$ consists of all vectors ${\bf x}=(x_1,\ldots,x_n) \in S$ for which there exist disjoint coordinate subsets $K$ and $L$ of $\{1,\ldots,n\}$ such that given the coordinates $x_i, i \in K$ one has that ${\bf x} \in A$ regardless of the values {of ${\bf x}$} on the remaining coordinates, and likewise that ${\bf x} \in B$ given the coordinates {$x_j, j \in L$}.
	For a finite product of discrete spaces endowed with a product measure, the BKR inequality
	\bea \label{eq:abs.BKR1}
	P(A \Box B) \le P(A)P(B)
	\ena 
	was conjectured by van den Berg and Kesten \cite{BK} and proved by Reimer \cite{Rei}.

	In \cite{GR} inequality \eqref{eq:abs.BKR1} was extended to general product probability spaces, replacing $A \Box B$ by 
	the set $A \eebbox B$ consisting of those outcomes {${\bf x}$} for which one can {only assure} with probability one that ${\bf x} \in A$ and ${\bf x} \in B$ based only on the revealed coordinates in $K$ and $L$ as above. A strengthening of the original BKR inequality \eqref{eq:abs.BKR1} results, due to the fact that $A \Box B \subseteq A \eebbox B$. In particular, it may be the case that  $A \Box B$ is empty,  while $A \eebbox B$ is not.

	We propose the further extension $A \stbbox B$ depending on probability thresholds $s$ and $t$, where $A \eebbox B$ is the special case where both $s$ and $t$ take the value one. The outcomes ${\bf x}$ in $A \stbbox B$ are those for which disjoint sets of  coordinates $K$ and $L$ exist such that given the values of $\bf x$ on the revealed set of coordinates  $K$, the probability that  $A$ occurs is at least $s$, and given the coordinates of $\bf x$ in $L$, the probability of $B$ is at least $t$.  We provide simple examples that illustrate the utility of these extensions.
\end{abstract}\bigskip

\section{Introduction}
\label{intro}
The `box' $A \Box B$ of events $A$ and $B$, consisting of the set of outcomes where $A$ and $B$ `occur for disjoint reasons,'  was first considered in \cite{vdBK} on the space $S = \{0,1\}^n$ with the uniform probability measure. More generally, on the product
$S=\prod_{i=1}^n S_i$ of general spaces $S_i,i=1,\ldots,n$,
with $[n]=\{1,\ldots,n\}$ and for $K \subseteq [n]$, define the $K$ cylinder of ${\bf x} \in S$ as 
\beas 
[{\bf x}]_K=\{{\bf y} \in S: y_i=x_i\,\, \forall i \in K\}.
\enas
The box operation is given by
\bea \label{def:classical_box}
A \Box B =\{{\bf x} \in S: \mbox{there exist disjoint $K,L \subseteq [n]$ such that $[{\bf x}]_K \subseteq A, [{\bf x}]_L \subseteq B$} \}.  
\ena
An equivalent definition of the box operation in \eqref{def:classical_box} in terms of 
cylinders $[A]_K$ of subsets $A$ of $S$, appearing as (2) in \cite{AGH} is
\bea \label{def:AGH.cyl.box1}
A \Box B = \bigcup_{0 \le |K| \le n} [A]_K \cap [B]_{K^c}\quad \mbox{where} \quad[A]_K=\{{\bf x }: [{\bf x}]_K \subseteq A\},
\ena
and $K^c$ denotes the complement of $K$.
The set $A \Box B$ represents the event that $A$ and $B$ occur for disjoint reasons, as for any ${\bf x} \in A \Box B$, the event $A$ can be verified to occur based only on the coordinates of ${\bf x}$ in $K$, regardless of the values of the remaining coordinates, and similarly so for coordinates in $L$ for $B$.  Clearly $A \Box B \subseteq A \cap B$. 
Note that if $K \subseteq M \subseteq [n]$ then $[{\bf x}]_M \subseteq [{\bf x}]_K$. This implies readily that in \eqref{def:classical_box}  we can assume without loss of generality that $L=K^c$. 
Note also that the choice of $K$ and $L$ in \eqref{def:classical_box} may depend on ${\bf x}$, and we may write $K({\bf x})$ and $L({\bf x})$ for specific choices. 

The famous {BKR} inequality, conjectured in \cite{vdBK} and proved in \cite{Rei} for the uniform probability measure on $\{0,1\}^n$, states that
\bea \label{eq:box.inequality}
P(A \Box B) \le P(A)P(B).
\ena
The validity of \eqref{eq:box.inequality} for the uniform measure on $\{0,1\}^n$ was shown to imply the inequality on a finite product of discrete probability spaces  in \cite{vdBF}.
{Henceforth we denote the end of an example by $\blacktriangle$.}

\begin{example} [A simple example: bond percolation.]\label{ex:simpleperc}  Consider the following graph:
	
	
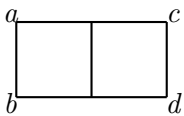
\centerline{\begin{tikzpicture}
		\draw [thick](0,0) --(1,0) -- (2,0);
		\draw [thick](0,1) --(1,1) -- (2,1);
		\draw [thick](0,0) --(0,1);
		\draw [thick](1,0) --(1,1);
		\draw [thick](2,0) --(2,1);
		\node at (-.08,1.08) {a};
		\node at (2.08,1.08) {c};
		\node at (-.08,-.08) {b};
		\node at (2.08,-.08) {d};
		\end{tikzpicture}}
	\noindent Edges are open or closed independently with probability $1/2$.
	Since there are 7 edges we may let $S=\{0,1\}^7$, and numbering the edges arbitrarily, the coordinates $x_i$ of ${\bf x}=(x_1,\ldots,x_7) \in S$ indicate whether edge $i$ is open or closed. Let $A=\{$there exists a path of open edges from $a$ to $c$$\}$.
	A direct count of the events that corresponds to $A$ shows that under the uniform distribution,  $P(A)=44/2^7$.  For $B=\{$there exists a path of open edges from $b$ to $d$$\}$ it is easy to see that $P(A \Box B)=2^3/2^7=1/16$  as $A \Box B$ is the set of all paths for which all bonds on the top and bottom rows of the graph are open.
	As $P(A)P(B)=({44/2^7})^2$,
	in this case `the excess multiple' in the bound \eqref{eq:box.inequality} is given by 
	\bea \label{eq:excess.multiple}
	\frac{P(A)P(B)}{P(A \Box B)}= \frac{({44}/{2^7})^2}{1/16} \approx 1.9.
	\ena 
	\qex
\end{example}

A version of the BKR inequality for a finite product of arbitrary probability spaces,  including finite spaces, and discrete or non-discrete infinite  spaces,  was considered in \cite{GR}. The
case where $A$ and $B$ are subsets of finite or countable products of ${\mathbb R}$, and more generally of Polish spaces, was considered in \cite{AGH}. This later work raises important issues 
regarding the measurability of $A \Box B$. In particular it is shown that the BKR combination $A \Box B$ of Lebesgue measurable sets need not be Lebesgue measurable. 
For an interesting practical application in a non discrete case see \cite{AGMS}.

We now turn to the consideration of $A \Box B$ in the framework of \cite{GR} and its slightly corrected version \cite{GRa}, on a finite product
$S=\prod_{i=1}^n S_i$ of general spaces $S_i,i=1,\ldots,n$ endowed with the product measure  $P=\prod_{i=1}^{n}P_i$ on the product sigma algebra of $S$ generated by given sigma algebras on $S_i$. 
The question of the measurability of $A \Box B$ did not arise in \cite{GR} since, as discussed in \cite{GRa}, the set $A \eebbox B$, which expresses the event that $A$ and $B$ occur for disjoint reasons almost surely, defined formally in \eqref{eq:specilize.function.version}, was considered instead.  The results of \cite{GR} show that the combination $A \eebbox B$ is measurable whenever $A$ and $B$ are, and that
\eqref{was.5} and \eqref{eq:oovsprod} below hold, that is 
\bea
\label{eq:allwekno}
A \Box B \subseteq A \eebbox B \quad \mbox{and} \quad P(A \eebbox B) \le P(A)P(B). 
\ena 
Further, in \cite{GR} it is shown that the latter inequality also holds when $P$ is the completion of any product measure. The proof in \cite{GR} of the inequality in \eqref{eq:allwekno} relies on Reimer's proof in \cite{Rei}. 
The implications of \eqref{eq:allwekno} are seen to be stronger than the BKR inequality \eqref{eq:box.inequality}. In particular, if $A \Box B$ is measurable, then it follows from \eqref{eq:allwekno} that $P(A \Box B) \le P(A)P(B)$.

Lemma 1 in \cite{AGH} shows that if $A$ and $B$ are Borel subsets of $[0, 1]^n$ then $A \Box B$ is Lebesgue measurable. As stated in Section 8 of \cite{AGH}, the same proof holds on the product of Polish subspaces of $\mathbb{R}$ equipped with the completion of any Borel product probability. Proposition \ref{prop:partof4} shows that  the qualifier `of $\mathbb{R}$' can be dropped. 

\begin{proposition}\label{prop:partof4}
	Let $S$ be a finite product of Polish spaces endowed with the product Borel sigma algebra and a product Borel probability measure $P$, and let $A$ and $B$ be Borel subsets of S.
	Then $A \Box B$ is measurable with respect to the completion $\overline{P}$ of $P$, and
	$\overline{P}(A \Box B) \le P(A)P(B)$. 
\end{proposition}

The measurability claim follows from definition \eqref{def:AGH.cyl.box1} and the result, shown in Section \ref{sec:meas}, the Appendix on measurability, that $[A]_K$ is  measurable with respect to the completed measure under the conditions of the Proposition; the inequality then follows by \eqref{eq:allwekno}. Here we have used common terminology where measurability of a set with respect  to the completion of $P$ means measurability with respect to the completion of the relevant sigma algebra with respect to $P$.

In Section \ref{sec:00box} we review the definition and meaning of the $\eebbox$ operation and provide some examples.
In Section \ref{sec:cyl} we look for an improvement of \eqref{eq:allwekno}, where for sets $A$ and $B$ we define in \eqref{def:overset2T}  minimal sets ${A}^{(\tiny{0})}$ and ${B}^{(\tiny{0})}$, `cylindrical cores',
satisfying ${A}^{(\tiny{0})} \subseteq A, \, {B}^{(\tiny{0})} \subseteq B$, and  $A \Box B = {A}^{(\tiny{0})} \Box {B}^{(\tiny{0})}$.
Measurability of  ${A}^{(\tiny{0})}$ under the conditions of Proposition \ref{prop:partof4} is shown in Proposition \ref{prop:LebegT}. The main purpose of this paper is achieved in Section \ref{sec:lenient}, which extends the set operation in \cite{GR} that captures the event that $A$ and $B$ occur for disjoint reasons almost surely, to one where $A$ and $B$ occur for disjoint reasons with at least some prescribed probabilities, and to provide a bound of BKR type. Examples are given in Section \ref{sec:lenexam}, some brief history on the BKR inequality in continuous spaces in Section \ref{sec:GuRa}, and various technical 
measurability issues are relegated to Section \ref{sec:meas}.

\section{The 11 Box Operation}\label{sec:00box}

With the notation given further meaning later, we {define} the box type operation $A \eebbox B$ on events $A$ and $B$. Apart from the notation $A \eebbox B$, these definitions appear in \cite{GR} and \cite{GRa}.  Unlike $A \Box B$, the operation $A \eebbox B$ depends on the measure $P$, taking into account probability zero events in a natural probabilistic way. The dependence on $P$ will be suppressed in the notation, here and in related definitions such as \eqref{def:stbox}.

Throughout this work, for each $i=1,\ldots,n$, let $(\spa_i,\bbs_i,P_i)$ be a probability
space, and set $\spa=\prod_{i=1}^n \spa_i$,
$\bbs=\bigotimes_{i=1}^n \bbs_i$,  $P=\prod_{i=1}^n P_i$ to be the product space, product sigma algebra, and product probability, respectively. For $K \subseteq [n]$ we set $S_K=\prod_{i \in K}^n \spa_i$, \, $\bbs_K=\bigotimes_{i \in K} \bbs_i$ and $P_K  =\prod_{i \in K} P_i$. For $A$ and $B$ in $\bbs$ define $A \Box B$ as in \eqref{def:classical_box} and note that it is shown in \cite{AGH} that the latter set may not be measurable. However, as a set it is well defined under the present general setting.

As in \cite{GR}, we consider real valued measurable functions on $\spa$, that is, functions 
that are $(\bbs,\mathbb{B})$ measurable, where $\mathbb{B}$
denotes the Borel sigma algebra of the real line ${\mathbb R}$. Let
$f$ and $g$ be two such functions, and for $K$ and 
$L$ subsets of $[n]$, let
\bea \label{eq:function.version}
\underline{f}_K({\bf x}) = {\rm ess} \inf_{{\bf y}
	\in [{\bf x}]_K}f({\bf y}), \quad \mbox{and} \quad
\underline{g}_L({\bf x}) = {\rm ess} \inf_{{\bf y} \in [{\bf
		x}]_L}g({\bf y}),
\ena
where the essential infimums defining $\underline{f}_K({\bf x})$ and $\underline{g}_L({\bf x})$ are taken with respect to the product probability measure on the coordinates in the complements $K^c$ and $L^c$ respectively.  
{See the Appendix of \cite{GR} where it is shown that  $\underline{f}$ and $\underline{g}$ of \eqref{eq:function.version}  are measurable, and also extensions for handling completions.}
Below we use $\max_{K \cap L = \emptyset}$ to indicate maximum over all pairs of disjoint sets $K,L \subseteq [n]$.

In Theorem \ref{thm1} below we consider random vectors ${\bf X}=(X_1,\ldots,X_n)$ where each $X_i$ is a measurable function from some probability space to the measurable
spaces $(\spa_i,\bbs_i)$, so that  $X_i \in S_i$, and for sets $B_i \in \bbs_i$ we have $P_i(B_i)=P(X_i \in B_i)$. We assume that  ${\bf X} = (X_1,\ldots,X_n)$ is distributed according to the product probability $P$, so that the components $X_i$ are independent.
Clearly, in moving to this framework no generality over that obtained by working on the basic probability space is lost. Moreover, function formulations such as in the theorem allow for statements that are more general than relations between the probabilities of sets; see \cite{GR} for some examples.

\begin{theorem}[{Goldstein and Rinott, 2007}]\label{thm1} 
	Let ${\bf X} = (X_1,\ldots,X_n) \in S$ be a random
	vector with a product distribution $P$. Then
	\bea \label{eq:finess} E\left\{\max_{K \cap L = \emptyset}\underline{f}_K({\bf X})\underline{g}_L({\bf
		X})\right\} \leq E\left\{f({\bf X})\right\}\,E\left\{g({\bf
		X})\right\}. \ena
\end{theorem}
We emphasize that Theorem \ref{thm1} applies to functions of a vector ${\bf X}$ having independent coordinates taking values in any measure
spaces.  

Specializing {\eqref{eq:function.version} and} \eqref{eq:finess} to the case where $f$ and $g$ are the indicator functions ${\bf 1}_A({\bf x})$ and ${\bf 1}_B({\bf x})$ of $A$ and $B$ respectively, 
we define $A \eebbox B$ by its indicator
\bea \label{eq:specilize.function.version}
{\bf 1}_{A \eebbox B}({\bf x})  = \max_{K \cap L = \emptyset}
\underline{{\bf 1}_A}_K({\bf x}) \underline{{\bf 1}_B}_L({\bf x}),\ena
Theorem \ref{thm1} now yields
\bea \label{eq:oovsprod} P(A \eebbox B) \le P(A)P(B).
\ena
We mentioned earlier that in \eqref{def:classical_box}  we can assume without loss of generality that $L=K^c$. We claim that this is true  also for \eqref{eq:specilize.function.version}, where again we can replace $K \cap L = \emptyset$ by $L=K^c$.  Clearly, such a replacement restricts the collection of sets over which the maximum is taken, so the quantity so obtained lower bounds the right hand side of \eqref{eq:specilize.function.version}. On the other hand, for any disjoint $K$ and $L$ we have $L \subseteq K^c$, which in turn implies $\underline{{\bf 1}_B}_L({\bf x}) \le \underline{{\bf 1}_B}_{K^c}({\bf x})$, as the infimum on the right is over a smaller set, see \eqref{eq:function.version}, implying
\beas
\underline{{\bf 1}_A}_K({\bf x}) \underline{{\bf 1}_B}_L({\bf x}) \le \underline{{\bf 1}_A}_K({\bf x}) \underline{{\bf 1}_B}_{K^c}({\bf x}).
\enas
In particular if the right-hand side of  \eqref{eq:specilize.function.version} takes the value 1, then for some disjoint $K$ and $L$ we have $\underline{{\bf 1}_A}_K({\bf x}) \underline{{\bf 1}_B}_L({\bf x})=1$ which implies $\underline{{\bf 1}_A}_K({\bf x}) \underline{{\bf 1}_B}_{K^c}({\bf x})=1$, so the replacement also upper bounds the right hand side of \eqref{eq:specilize.function.version}, and our claim follows.

We now consider definition \eqref{def:classical_box} of $A \Box B$ in general spaces, which, as noted above, is a set operation that does not depend on the underlying measure, in contrast to $A \eebbox B$.
As shown in   (5) of \cite{GRa}, we have the inclusion
\bea \label{was.5}
{A \Box B \subseteq A \eebbox B}.
\ena
Now \eqref{eq:oovsprod} and \eqref{was.5} imply  \eqref{eq:box.inequality}  in the general setting  above provided $A \Box B$ is measurable. To see \eqref{was.5}, note that replacing essential infimum by infimum in \eqref{eq:function.version},  inequality \eqref{was.5} becomes equality. Hence \eqref{was.5} holds as stated because  the essential infimum is at least as large as the infimum. In words, elements of $A \Box B$
demand that $A$ and $B$  hold for all outcomes satisfying disjoint reasons for their occurrence, while $A \eebbox B$ requires only that the `reasons' imply that $A$ and $B$ hold almost surely.

To develop some intuition on the $\eebbox$ box type operation, consider two individuals, one of whom desires that event $A$ happens, while the other desires $B$; we label the individuals according to their desired event with little danger of confusion. If ${\bf x} \in A \Box B$ then disjoint sets of coordinates of ${\bf x}$ may be revealed to $A$ and $B$ so that each will know with certainty that their event has occurred. Definition \eqref{eq:specilize.function.version} extends the box operation to accommodate situations where the disjoint sets of coordinates of $\bf x$ revealed to $A$ and $B$ indicate only that their event has occurred almost surely, that is, with probability one, but which may not certify it absolutely. It is easy to see that unlike $A \Box B$, the set $A \eebbox B$ is not necessarily contained in $A \cap B$; an example that follows gives one such instance.

It is natural to ask what might be the differences, and possible advantages, of defining the box operation $\eebbox$ as in \eqref{eq:specilize.function.version}, and more generally for functions as in \eqref{eq:function.version} as {considered} in Theorem \ref{thm1},  where the underlying measure is involved, in contrast to the more `straightforward' extension that preserves \eqref{def:classical_box} as a pure set operation.
For one, the approach taken in Theorem \ref{thm1} avoids certain subtle measurability difficulties that arise in the `straightforward' approach. In particular, Example 2 in \cite{AGH} presents a situation where the classical box of two Lebesgue measurable sets fails to be Lebesgue measurable. 
In contrast, {as shown in \cite{GR}}, the function operation as defined in \eqref{eq:function.version} preserves measurability.

Technical matters aside, the approach of Theorem \ref{thm1} also has advantages from a probabilistic perspective, as illustrated by the following simple example, somewhat related to Example 2 in \cite{AGH}.
Let 
\bea \label{eq:soda1}
S = [0,1]^2 \qmq{and}
A=B = \{(x_1,x_2) \in [0,1]^2: x_1 \not = x_2\}.
\ena
It is easy to see that $A \Box B$ as defined in \eqref{def:classical_box} is empty. For instance, for any ${{\bf x}} \in [0,1]^2$, the cylinder $[{{\bf x}}]_{\{1\}}$ is not contained in $A$ (only) because it contains the point $(x_1,x_1)$, which is not in $A$. However, for any continuous measure on $S$, this single point is a set of measure zero, with the result that $A \eebbox B = [0,1]^2$.  

Here is an intuitive way to think about this example: suppose Alice and Bob each need to choose a single instant of time in the interval $[0,1]$ required for completion of a certain task. Assume there is only a single resource, and both will fail if and only if they require the resource at the same precise instant. To succeed for disjoint reasons as in \eqref{def:classical_box} means that each one can be sure of their own success independently of the other's choice of time. That here $A \Box B = \emptyset$ is a reflection of the fact that this is impossible. However, if they select times independently
by any continuous distribution
in the unit interval, then each one of them will succeed with probability one for any choice of the other. And indeed, for this case $A \eebbox B = [0,1]^2$. This example, though simple, illustrates how sets of measure zero, which typically are a technical nuisance that can effectively be ignored, can greatly affect box type operations, here in particular causing the maximum possible discrepancy in the inequality
\beas
0=P(\emptyset)=P(A \Box B) < P(A \eebbox B)=1.
\enas

\begin{example}[Continuous percolation]\label{ex:perc}
	Consider the random geometric graph ${\cal G}$ whose vertex set is a
	collection of $n$ independent points ${{\bf X}}={\{}X_1, \ldots, X_n{\}}$, each sampled from a distribution having a strictly positive density  on $[0,1]^2$. An edge is drawn between distinct points $X_i$ and $X_j$ if $||X_i-X_j||\le 2r$, where $r>0$ and $|| \cdot ||$ is any norm in $\mathbb{R}^2$. An equivalent way to view the graph is to connect two distinct points when the two circles of radius $r$ having these points as their centers intersect. 
	For percolation models based on randomly chosen balls see, for example, \cite{MP} and \cite{MR} and \cite{LPZ}.

	Let $A$ be the event that there exists a `path along edges from left to right of $[0,1]^2$, that is, a path starting from a point that is within distance of at most $r$ from the $y$-axis, and ending at a point within distance at most $r$ from the line $x=1$. Similarly, let $B$ denote the event that there exists a path from  bottom to top. 
	It is easy to see that for $n$ large, with positive probability, there exist realizations that allow crossing from left to right and from bottom to top on disjoint paths, that is, such that the two paths do not share a vertex.
	In this case the event that $A$ and $B$ occur for disjoint reasons, which coincides with the event $A \Box B$, is non-empty.  These events can be seen as continuous versions of those studied in \cite{BK}.\qex
\end{example}

\begin{example}[Continuous percolation with annihilation]\label{ex:percann}
	We next introduce annihilation to the percolation model on the graph ${\cal G}$ of Example \ref{ex:perc}, and return to this setting again in Section \ref{sec:lenexam}.
	If for some $i \not =j$ the  points $X_i$ and $X_j$ coincide, that is, are equal, then any edge attached to these points is annihilated. The coinciding points themselves are not annihilated, so the number of points remains $n$.   Processes with annihilation are ubiquitous, appearing in connections to random walks, Brownian motion, branching processes and more, based on motivations from physics and chemistry. As $X_i,i=1,\ldots,n$ are independent, each with a continuous distribution, the probability of annihilation is zero. Let $A$ and $B$ be as defined in Example \ref{ex:perc}. The fact that $A \Box B =  \emptyset$ in the presence of annihilation, hence making inequality \eqref{eq:box.inequality} trivial, is covered by the following more general result in the context of this model.$\qex$
\end{example}

\begin{proposition}\label{prop:1} 
	For the percolation model  with annihilation $A \Box B=\emptyset$, where $A$ and $B$ are any events that require that the edges of at least one point are not annihilated.
\end{proposition} 
\proof For any nontrivial event $A$ that requires that edges of at least one point  not be annihilated, any ${{\bf x}}$ and any index set $K$ such that $|{K^c}| \ge |K|$, it is impossible that $[{{\bf x}}]_K \subseteq A$ as the {$|{K^c}|$} unspecified coordinates in $[{{\bf x}}]_K$, those indexed by the larger set ${K^c}$, can annihilate all the edges attached to the points indexed by $K$. Hence  $A \Box B$ must be empty, as given complementary sets $K$ and $L$, at least one set, say $K$, satisfies $|K^c| \ge |K|$. \qed 

We remark that Proposition \ref{prop:1} has a flavor similar to the example based on \eqref{eq:soda1}, in that the outcome of the classical box operation may be determined by events that occur with probability zero. Returning to the  sets $A$ and $B$ in Example \ref{ex:percann}, Proposition \ref{prop:1} implies $A \Box B=\emptyset$. However, for $n$ sufficiently large, $P(A \eebbox B)>0$, and moreover, this probability converges to one as $n \rightarrow \infty$. To see this, consider a finite set of circles $C_i$, $i=1,\ldots,m$ of radius $r/2$ and centers in $[0,1]^2$ such that $C_1$ intersects the $y$-axis, and $C_m$ intersects the line $x=1$, and $C_{i} \cap C_{i+1} \ne \emptyset$ for $i=1,\ldots,m-1$. As $n \rightarrow \infty$, for all $i=1,2,\ldots,m$ there will be a point $X_i \in C_i$ with probability approaching one, and the edges between these points will not be annihilated with probability one. Hence, these points will form a path from left to right corresponding to the event $A$, and a similar argument for $B$ implies that $P(A \eebbox B) \rightarrow 1$.

Finally we comment that  ${\bf x} \in A \Box B$ under the model without annihilation if and only if ${\bf x} \in A \eebbox B$ under the model with annihilation, that is, when edges attached to coincident points are annihilated. In other words,  $A \Box B=A \eebbox B$, where the  set on the left-hand side is defined for the model without annihilation, and the right-hand set is defined for the model with annihilation.

\section{Cylindrification }\label{sec:cyl}
Before proceeding to our extension of the box operation, we pause to make a simple observation about potential improvements to the original box inequality, and our subsequent generalizations. Here we show that `cylindrical cores' of the sets that comprise the box operation may be separately extracted, yielding no loss to the left-hand side of \eqref{eq:box.inequality}, but a possibly smaller quantity on the right. For this purpose, the cylindrical core $E^{(0)}$ of a given $E \in \bbs$, that extracts the non-trivial cylinders that $E$ contains, is defined for $n \ge 1$ by
\bea \label{def:overset2T}
{E}^{(\tiny{0})}=\bigcup_{{({\bf x},K)}\,:\,{K \subseteq [n]},0 \le |K| \le n-1, [{\bf x}]_K \subseteq E} [{\bf x}]_K.
\ena
Note that for any ${\bf x} \in S$ we have  ${\bf x} \in [{\bf x}]_K \subseteq S$ for any $K$, and therefore ${\bf x} \in {S}^{(\tiny{0})}$,
implying that   ${S}^{(\tiny{0})}=S$.
Recalling the definition of $[A]_K$ in \eqref{def:AGH.cyl.box1}, it is easy to verify that \eqref{def:overset2T}  is equivalent to
\bea \label{def:overset3T}
{E}^{(\tiny{0})}=\bigcup_{{K \subseteq [n]},0 \le |K| \le n-1} [E]_K.
\ena

In Proposition \ref{prop:coreT} we show that ${E}^{(\tiny{0})}$ is the smallest set that preserves the result of the box operation applied to $E$ and any proper subset of $S$.

Since for all $E$ it holds that $E \Box S = E$, the equality $A \Box S = A^{(\tiny{0})} \Box S$ would imply $A=A^{(\tiny{0})}$. Thus the claim \eqref{eq:AboxB.overset.AT} has been restricted to proper subsets $B$ of $S$.
Recalling $\spa=\prod_{i=1}^n \spa_i$, we assume for the last part of the next proposition that for all $i=1,\ldots,n$ the set $\spa_i$ is non-trivial in that it contains more than a single point.
\begin{proposition}\label{prop:coreT}
	{For all $A \subseteq S$} the set ${A}^{(\tiny{0})}$ defined in \eqref{def:overset2T} satisfies ${A}^{(\tiny{0})} \subseteq A$, and 
	\bea \label{eq:AboxB.overset.AT} 
	A \Box B = A^{(\tiny{0})} \Box B \quad \forall B \not = S.
	\ena 
	If $A$ and $B$ are proper subsets of $S$, then
	\bea \label{eq:AboxBsuoo}
	A \Box B = A^{({\tiny{0}})} \Box B=A \Box B^{({\tiny{0}})} = A^{({\tiny{0}})} \Box B^{({\tiny{0}})}.
	\ena
	If the space $\spa=\prod_{i=1}^n \spa_i$ satisfies $|S_i| >1$ for all $i=1,\ldots,n$ then ${A}^{(\tiny{0})}$  is the smallest set that satisfies
	\eqref{eq:AboxB.overset.AT} in the sense that   if \eqref{eq:AboxB.overset.AT} also holds for a set ${A}^{(\tiny{1})}$ in place of $A^{(\tiny{0})}$ then $A^{(\tiny{1})} \supseteq A^{(\tiny{0})}$.
\end{proposition}

\noindent {\em Proof: } {First it is clear that $A^{(\tiny{0})}$, being a union of subsets of $A$, is a subset of $A$.}
Hence, by the evident monotonicity of the box operation, in order to prove \eqref{eq:AboxB.overset.AT} we are required only  to show that
$A \Box B  \subseteq {A}^{(\tiny{0})} \Box B$ for all $B \not = S$.
If ${\bf x} \in A \Box B$ then there exist disjoint $K$ and $L$ such that $[{\bf x}]_K \subseteq A$ and $[{\bf x}]_L \subseteq B$. Note that $B \not = S$ implies $|L| \ge 1$ and hence $|K| \le n-1$. Now by definition $[{\bf x}]_K \subseteq A^{(\tiny{0})}$, and ${\bf x} \in A^{(\tiny{0})} \Box B$.

The first and third equalities in \eqref{eq:AboxBsuoo} follow by using \eqref{eq:AboxB.overset.AT} and that $B$ is a proper subset of $S$, that $A^{({\tiny{0}})} \subseteq A$ and that $A$ is a proper subset of $S$.  The second equality then follows by reversing the roles of $A$ and $B$.

We finally show the minimality of $A^{(\tiny{0})}$ when $|S_i|>1$ for all $i=1,\ldots,n$. If ${\bf x} \in A^{(\tiny{0})}$ then for some $K$ satisfying $|K|\le n-1$ we have $[{\bf x}]_K \subseteq A^{(\tiny{0})}$. Now, with $L={K^c}$ let $B = [{\bf x}]_L$; then ${\bf x} \in A^{(\tiny{0})} \Box B$. Since $|L|=n-|K| \ge 1$ and $|S_i|>1$ {for all $i=1,\ldots,n$}  we have that $B \not = S$ and therefore ${\bf x} \in A^{(\tiny{0})} \Box B =  A^{(\tiny{1})} \Box B$, so in particular ${\bf x} \in A^{(\tiny{1})}$. The inclusion $ A^{(\tiny{0})} \subseteq  A^{(\tiny{1})}$ is hence demonstrated.  
\qed 
Applying this result to any proper subsets $A$ and $B$ of $S$, we obtain, for possible improvements in \eqref{eq:box.inequality},
\bea \label{eq:AboxB.oversetT}
{A}^{(\tiny{0})} \subseteq A, \quad {B}^{(\tiny{0})} \subseteq B, \lcolor{\quad}  A \Box B = {A}^{(\tiny{0})} \Box {B}^{(\tiny{0})} \qmq{and} P(A\Box B)\le  P(A^{({\tiny{0}})}) P(B^{({\tiny{0}})}),
\ena
with the last inequality holding provided all sets involved are measurable;  see Proposition \ref{prop:LebegT}.

For the rest of this section we focus on the case where $\spa=\prod_{i=1}^n \spa_i$, a product of Polish spaces, $\bbs$ is the Borel product sigma algebra on $\spa$, $P$ is a product Borel measures on $S$ and the sets $A$ and $B$ considered below are Borel.  
Recall that Proposition \ref{prop:partof4} states  that  $A \Box B$ is measurable with respect to the completion $\overline{P}$ of $P$  in this framework.  Using similar arguments, the next proposition gives the measurability of $A^{({\tiny{0}})}$ under the same conditions.

\begin{proposition}\label{prop:LebegT}
	Let $S$ be a finite product of Polish spaces endowed with the product Borel sigma algebra and a product Borel probability measure $P$, and let $A$ and $B$ be Borel subsets of S. Then the sets 
	$A^{({\tiny{0}})}$ and 
	$B^{({\tiny{0}})}$
	are measurable with respect to the completion $\overline{P}$ of $P$, and  
	for all $A$ and $B$ proper Borel subsets of $S$, we have
	$$\overline{P}(A \Box B) \le {\overline P}(A^{({\tiny{0}})}) {\overline P}(B^{({\tiny{0}})}).$$
\end{proposition}
\noindent {\em Proof: }
We prove that the set $[A]_K$, defined in \eqref{def:AGH.cyl.box1},  is $\overline{P}$ measurable in Section \ref{sec:meas}, the Appendix on measurability. The measurability claim now follows from \eqref{def:overset3T}.

To prove the inequality,  by Proposition \ref{prop:coreT} and \eqref{was.5} we have
$A \Box B {=} A^{({\tiny{0}})} \Box B^{({\tiny{0}})}\subseteq A^{({\tiny{0}})} \eebbox B^{({\tiny{0}})}$, and the result now follows by \eqref{eq:oovsprod}. 
\qed

We now illustrate the $E \mapsto E^{(\tiny{0})}$ operation with the following simple example.
\begin{example}[Bond Percolation]
	Consider bond percolation on the graph below. There are 13 edges, each of which is open or closed independently, yielding $|S|=2^{13}$. A path from left to right in this example means a path from either $a$ or $b$ to either $c$ or $d$ along open edges. 
	
	
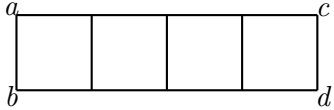
\centerline{\begin{tikzpicture}
		\draw [thick](0,0) --(1,0) -- (2,0) -- (3,0) -- (4,0);
		\draw [thick](0,1) --(1,1) -- (2,1) -- (3,1) -- (4,1);
		\draw [thick](0,0) --(0,1);
		\draw [thick](1,0) --(1,1);
		\draw [thick](2,0) --(2,1);
		\draw [thick](3,0) --(3,1);
		\draw [thick](4,0) --(4,1);
		\node at (-.08,1.08) {a};
		\node at (-.07,-.07) {b};
		\node at (4.07,-.07) {d};
		\node at (4.07,1.07) {c};
		\end{tikzpicture}}
	Let $A=$$\{$there exists a path from left to right$\}\cup$$\{$all edges are closed$\}\cup$$\{${exactly three edges are open}$\}$.
	In this case $A^{(\tiny{0})}=\{$there exists a path from left to right$\}$ which {is a strict subset} of $A$. More specifically, $A \setminus A^{(\tiny{0})}$ contains $1+ {13 \choose 3}$ points. Assuming, for example, a uniform probability, we obtain the gap $P(A)-P(A^{(\tiny{0})})=(1+ {13 \choose 3})/2^{13}=0.035$. It is easy to increase the gap by taking the union of $A$ with the set $\{$exactly two edges are open$\}$, or that exactly a particular set of edges is open, for any set that does not guarantee a path from left to right. \qex
	
\end{example}

\section{The Lenient $s,t$ Box} \label{sec:lenient}

In Section \ref{intro} we considered the box operation {$\eebbox$} that is not influenced by events of zero probability. Here, for $(s,t) \in [0,1]^2$,  we define the more general $\stbbox$ operation from which {$\eebbox$} can be recovered as the 
special case where $(s,t)={(1,1)}$. The parameters $s$ and $t$  give thresholds for the   `leniency' on certain conditional probabilities that $A$ and $B$ occur.
Our notations and definition of conditional probability  follow \cite{Brei} Chapter 4, and will be explained in detail below.

We now describe the `lenient' set operation $A \stbbox B$. Consider a random vector ${\bf X} = (X_1,\ldots,X_n)$ $\in S$ under the conditions of Theorem \ref{thm1}, and the probability that ${\bf X}$ is in $A$ conditioned on the values of the coordinates of ${\bf X}$ in some $K \subseteq [n]$, expressed as $P({\bf X} \in {A} \mid {\bf X}_{K}={\bf x}_{K})$, where for $K \subseteq [n]$ and for ${\bf x} \in S$ the vector ${\bf x}_K$ is the projection of ${\bf x}$ on $K$, that is, the $|K|$-vector consisting of the coordinates of $\bf x$ in $K$.  If this conditional probability is at least $s$ then revealing that the coordinates in $K$ of ${\bf X}$ and ${\bf x}$ agree implies that the conditional probability that ${\bf X} \in A$ is at least $s$.  
We define 
\begin{multline}\label{def:firstst}
	$A\stbbox B$=\{{\bf x} \in S: \exists K,L \subseteq [n],\,\,
	K \cap L = \emptyset \,\, \mbox{such that} \,\,\, P({\bf X} \in A \mid {\bf X}_{K}={\bf x}_{K})\ge s\,\,\, \mbox{and} \,\, P({\bf X} \in B \mid {\bf X}_{L}={\bf x}_{L})\ge t \}.
\end{multline} 
In particular, when both $s$ and $t$ are close to one, events $A$ and $B$ are both conditionally very likely to occur for disjoint reasons.

{The set $A \stbbox B$ as described is defined using conditioning on ${\bf X}={\bf x}$, rather than on a sigma algebra; we refer the reader to Definition 4.18. of \cite{Brei}. We also note, by Proposition 4.9 of \cite{Brei}, that when the conditioning  sigma algebra is generated by ${\bf X}$,  conditional expectation may be written as a function of ${\bf X}$. See also Corollary 4.38 of \cite{Brei}, which is also relevant for our case. Measurability of the sets $E_{r,K},E_r$ and $A \stbbox B$ introduced in this section are handled in Proposition \ref{prop:measure.st.sets}.

	Our formal definition of the $\stbbox$ operation requires specifying a particular version of the conditional probabilities in \eqref{def:firstst}.  For $K \subseteq [n]$, ${\bf y} \in S_K$ and ${\bf v} \in S_{K^c}$ let
	$\langle {\bf y},{\bf v}\rangle_K$ be the point in $S$ with $i^{th}$ coordinate given by	
	$$\left({\langle {\bf {y}},{\bf v} \rangle_K} \right)_i = \left\{
	\begin{array}{cc}
	{y}_i &  i \in K\\
	v_i & i \in {K^c},
	\end{array}
	\right.$$
	and for $E \in \mathbb S, 0 \le r \le 1$ and 
	$K {\subseteq} [n]$ let
	\bea \label{def:uv.cond}
	E_{r,K}=\left\{
	\begin{array}{cl}
		\left\{{\bf x} \in S: P_{K^c}\left(\{{\bf v} \in S_{K^c}: {\langle {\bf x}_{{K}}, {\bf v}  \rangle_K} \in E\}\right) \ge r \right\}  & \qm{for $|K| \le n-1, 0 \le r \le 1$} \\
		E & \qm{for $K=[n], 0<r \le 1$, and}\\
		S & \qm{for $K=[n], r = 0$.}
	\end{array}
	\right.
	\ena

	We note that the relation ${\bf x} \in E_{r,K}$ is equivalent to $[{\bf x}]_K \in E_{r,K}$ since it depends on ${\bf x} \in S$ only through ${\bf x}_K$, as the probability appearing in its definition only so depends. To describe the set $E_{r,K}$ in words, given ${\bf x} \in S$,  consider the set of  vectors ${\bf v} \in S_{K^c}$ such that when substituting  the coordinates of ${\bf x}$ with indices in ${K^c}$ by those of ${\bf v}$ yields a vector that is in $E$. If the ${K^c}$ marginal probability of the set of such ${\bf v}$ is at least $r$, then  ${\bf x}$ is in $E_{r,K}$.
	Illustrative special cases of the application of definition \eqref{def:uv.cond} are $E_{0,K}=S$ for all $K$ as the substitution of any coordinates in $K^c$ will yield a non-negative conditional probability. Also $E_{r,\emptyset} = S$ if $P(E) \ge r$, as even having no coordinates revealed one is still assured that $E$ will occur with probability at least $r$,
	and otherwise $E_{r,\emptyset} = \emptyset$, as $\langle {\bf x}_\emptyset, {\bf v}\rangle_\emptyset ={\bf v}$, which is in $E$ if and only if ${\bf v} \in E$, whose probability does not reach the threshold $r$.

	With $A_{s,K}$ and  $B_{t,L}$ defined as in \eqref{def:uv.cond} for measurable set $A$ and $B$ of $S$, and $(s,t) {\in} [0,1]^2$, we have
	\bea \label{def:stbox}
	A \stbbox B = \bigcup_{K \cap L = \emptyset} {A_{s,K} \cap B_{t,L}}.
	\ena
	
	Definitions \eqref{def:firstst} and \eqref{def:stbox} are equivalent.
	This equivalence becomes clear when noting, as stated in Proposition \ref{prop:measure.st.sets},
	that $P_{K^c}\left(\{{\bf v} \in S_{K^c}: \langle {\bf {y}}, {\bf v}  \rangle_K {\in} E\}\right)$,  defined for all ${\bf y}\in S_K$, is a version of the conditional probability $P({\bf X} \in E|{\bf X}_K={\bf {y}})$. 
	The event $A \stbbox B$ consists of all outcomes ${\bf x}$ for
	which there are disjoint sets of coordinate indices $K$ and $L$ such that the conditional probabilities of $A$ given the values of  $x_i$ for $i \in K$, and $B$  given the values of  $x_j$ for $j \in L$, are at least $s$ and $t$, respectively.

	It is  easy to see that for any events $A$ and $B$ we have $A \oobbox B=S$.  Moreover, 
	if $P(A)\ge s$  and   $P(B)\ge t$  then  $A \stbbox B = S$, as can be seen by choosing $K=L=\emptyset$. However, this choice is impossible if we impose $L=K^c$, and in fact, it is easy to see that the relations $P(A)\ge s$  and   $P(B)\ge t$ do not imply $A \stbbox B = S$ in general if we impose $L=K^c$ in \eqref{def:stbox}. This shows that
	in contrast to definition \eqref{def:classical_box} of $A \Box B$ and definition \eqref{eq:specilize.function.version} of $A \eebbox B$, in general $A \stbbox B$ of definition \eqref{def:stbox}   is not equivalent  to the one obtained by replacing the condition that $K \cap L = \emptyset$ by $L={K^c}$.
	See Example \ref{ex:threeside} below for another illustration. Nevertheless, in Example \ref{ex:coin} below, assuming that $K$ and $L$ are complementary in the definition of the set $A \stbbox B$ will not affect the set defined.

	The next proposition shows that definitions \eqref{def:stbox} and  \eqref{eq:specilize.function.version} are consistent, that is, that they coincide when $(s,t)=(1,1)$, thus explaining our earlier choice of notation for $A\eebbox B$.
	\begin{proposition}\label{prop:easyoo} The set $A \stbbox B$ obtained by specializing definition \eqref{def:stbox} for $(s,t)=(1,1)$ agrees with the set $A \eebbox B$ whose indicator is given in \eqref{eq:specilize.function.version}.
	\end{proposition}
	\proof By  \eqref{eq:specilize.function.version},
	${\bf x} \in A\eebbox B$ if and only if for some disjoint $K$ and $L$,
	${\rm ess} \inf_{{\bf y}	\in [{\bf x}]_K}{\bf 1}_A({\bf y})={\rm ess} \inf_{{\bf y}	\in [{\bf x}]_L}{\bf 1}_B({\bf y})=1$, 
	if and only if for some disjoint $K$ and $L$, $P_{K^c}\left(\{{\bf v} \in S_{K^c}: \langle {\bf x}_{{K}}, {\bf v}  \rangle_K \in A\}\right)=1$ and \\ $P_{L^c}\left(\{{\bf v} \in S_{{L}^c}: \langle {\bf x}_{{L}}, {\bf v}  \rangle_{{L}} \in B\}\right)=1$. These latter two equalities hold if and only if ${\bf x} \in A\stbbox B$ for ${(s,t)=(1,1)}$.  \bbox

	Theorem \ref{thm:bound.by.enlarged} below gives a version of the BKR inequality for the lenient box operation.  Given $E \in \mathbb S$ and $r {\ge} 0$, define the $r$-inflated set
	\begin{equation}\label{eq:nowinflate}
		E_r= \bigcup_{K \subseteq [n]} E_{r,K}.
	\end{equation}
	Since $E_{r,[n]}=E$ for $r>0$ and equals $S$ otherwise, it follows that $E \subseteq E_r$. Below we denote the  set $(E_r)^{(\tiny{0})}$ by $E_{{r}}^{(\tiny{0})}$.
	As $E_r$ is measurable for a measurable $E$  by Proposition \ref{prop:measure.st.sets} below, so are the sets $A_s$ and $B_t$ appearing in Theorem \ref{thm:bound.by.enlarged}, and  by Proposition \ref {prop:LebegT}, the sets $A_s^{(\tiny{0})}$ and $B_t^{(\tiny{0})}$ are measurable with respect to the completion of $P$ when $S$ is a finite product of Polish spaces endowed with the product Borel sigma algebra, on which $P$ is a product probability measure.

	\begin{theorem} \label{thm:bound.by.enlarged}
		Let $A$ and $B$ be in $\mathbb S$  and $(s,t) \in  [0,1]^2$.  Then with 
		$A_s$ and $B_t$
		as in \eqref{eq:nowinflate},
		\bea \label{eq:BKRst}
		P(A \stbbox B)   \le P(A_s)P(B_t).
		\ena
		If the sets $A_s^{(\tiny{0})}$ and $B_t^{(\tiny{0})}$ are measurable then 
		\bea \label{eq:nowinflate0}
		P(A \stbbox B)  \le P(A_s^{(\tiny{0})})P(B_s^{(\tiny{0})}) \le P(A_s)P(B_t).
		\ena
	\end{theorem}

	\noindent {\em Proof:}   We first show that
	$A \stbbox B \subseteq A_s \Box B_t $. Indeed, if ${\bf x} \in A \stbbox B$ then there exist disjoint $K$ and $L$ such that 
	\beas
	P_{K^c}\left(\{{\bf v} \in S_{K^c}: \langle {\bf x}_{{K}}, {\bf v}  \rangle_K \in A\}\right) \ge s \qmq{and} P_{L^c}\left(\{{\bf v} \in S_{L^c}: \langle {\bf x}_{{L}}, {\bf v}  \rangle_L \in B\}\right) \ge t.
	\enas
	In particular, $[{\bf x}]_K \subseteq A_{s,K} \subseteq A_s$ and $[{\bf x}]_L  \subseteq B_{t,L} \subseteq B_t$, and as $K$ and $L$ are disjoint we have ${\bf x} \in A_s \Box B_t$ as desired. Together with \eqref{was.5} we now have $A \stbbox B \subseteq A_s \Box B_t \subseteq A_s \eebbox B_t $, and now \eqref{eq:oovsprod} implies \eqref{eq:BKRst}. 
	
	To show \eqref{eq:nowinflate0}, we have 
	$$A \stbbox B \subseteq A_s \Box B_t= A_s^{(\tiny{0})} \Box B_t^{(\tiny{0})} \subseteq A_s^{(\tiny{0})} \eebbox B_t^{(\tiny{0})}$$
	where we first apply the inclusion just shown, followed by
	\eqref{eq:AboxB.oversetT} and  \eqref{was.5}. We now obtain the first inequality in \eqref{eq:nowinflate0} applying \eqref{eq:oovsprod}; the second one follows from the first two relations in \eqref{eq:AboxB.oversetT}.  \qed

	Measurability issues of the sets appearing in this section are summarized in the following proposition, which is proved in Section \ref{sec:meas},  the Appendix on measurability. Below we set ${\bf y}={\bf x}_{K} \in S_K$, and the set $\{{\bf v} \in S_{K^c}: {\langle {\bf x}_{{K}}, {\bf v}  \rangle_K} \in E\}$ becomes $\{{\bf v} \in S_{K^c}: {\langle {\bf y}, {\bf v}  \rangle_K} \in E\}$.
	\begin{proposition} \label{prop:measure.st.sets}
		The set $\{{\bf v} \in S_{K^c}: {\langle {\bf y}, {\bf v}  \rangle_K} \in E\}$, whose $P_{K^c}$ probability appears in \eqref{def:uv.cond}, is measurable, and $P_{K^c}\left(\{{\bf v} \in S_{K^c}: {\langle {\bf y}, {\bf v}  \rangle_K} \in E\}\right)$ defined for all ${\bf y}\in S_K$,  is a version of the conditional probability $P({\bf X} \in E|{\bf X}_K={\bf {y}})$. For
		a measurable set $E$, the sets $E_{r,K}$ of \eqref{def:uv.cond} and $E_r$ of \eqref{eq:nowinflate} are measurable. Also, for measurable sets $A,B$ the set $A \stbbox B$ of \eqref{def:stbox}  is measurable. Furthermore, the replacement of our version of the conditional probability by another in \eqref{def:uv.cond} changes  $E_{r,K}$ by a $P$-null set,  and hence the same applies to $A \stbbox B$, and  our results hold independently of the version chosen.  
	\end{proposition}

	\section{Lenient Box Examples} \label{sec:lenexam}
	\begin{example}\label{ex:coin}[An odd coin tossing problem]\label{example:odd.coin}
		A fair coin is to be tossed independently $n=2m+1$ times  for some integer $m \ge 1$. Let $A$ be the event that the first $m+1$ tosses are all heads, and $B$ the event that the last $m+1$ tosses are all tails. With $s=t=1/2$ it is not hard to see that $A \stbbox B$ consists of the following sequences: (a) the two sequences in which the first $m$ tosses are heads, the last $m$ tosses are tails and the middle toss could be either; (b) the $m$ sequences in which the $m+1^{st}$ toss is a head, the first $m$ tosses are all heads except for exactly one, and the last $m$ tosses are all tails, and (c) the $m$ sequences in which the $m+1^{st}$ toss is a tail, the first $m$ tosses are all heads, and the last $m$ tosses are all tails except for exactly one. Altogether there are $2(m+1)$ such sequences, hence
		\beas
		P(A \stbbox B) = 2 (m+1) (1/2)^n=(m+1)(1/2)^{2m}. 
		\enas 
		To illustrate the notions defined in previous sections, we determine the sets $A_{1/2,K}$ and $B_{1/2,L}$ as defined in \eqref{def:uv.cond}, and compute $P(A_{1/2})$.  For $A_{1/2,K}$ to contain some  $\bf x$, its projection  ${\bf x}_K$ must either specify all the first $m+1$ coordinates as heads, so that $A$ occurs with probability one, or heads in all but one of these coordinates, so that the probability that the unspecified coordinate is heads, and therefore that $A$ occurs, is $1/2$.  In particular, $A_{1/2,K} = \emptyset$ if $|K| < m$, so contributions to $A_{1/2}$ can only arise from sets $K$ satisfying $|K| \ge m$.
		Hence, $A_{1/2}$, following definition \eqref{eq:nowinflate},
		consists of $A$ and all sequences that have at most one tail in the first $m+1$ positions, and thus we obtain
		\beas
		P(A_{1/2})=(1/2)^{m+1}+ (m+1)(1/2)^{m+1} = (m+2)(1/2)^{m+1}.
		\enas 
		
		Arguing similarly for $B$, the inequality 
		\beas
		(m+1)(1/2)^{2m} = P(A \stbbox B)  \le P(A_{1/2})P(B_{1/2})= \{(m+2)(1/2)^{m+1}\}^2=(m+2)^2 (1/2)^{2m+2}
		\enas
		is consistent with  the conclusion of Theorem \ref{thm:bound.by.enlarged}. These calculations show that in this case the excess multiple in the bound, such as the one computed in \eqref{eq:excess.multiple} for the standard box, is given by
		\beas
		\frac{P(A_{1/2})P(B_{1/2})}{P(A \stbbox B)} = \frac{(m+2)^2}{4(m+1)}.
		\enas

		In this case the lenient box is not empty, yet not only is the box empty, but also the intersection, that is 
		$A \cap B = \emptyset$. Intuitively this means that
		$A$ and $B$ can never occur together and certainly not for disjoint reasons.  Thus if Alice is betting on $A$ and Bob on $B$, then it is impossible for both of them to win their bets. However, certain revealed coordinates can make each of them hopeful. The set $A \stbbox B$ with {$(s,t)=(1/2,1/2)$} contains those ${\bf x}$'s for which there exists a set of coordinates such that if Alice conditions on them her probability of winning the bet is at least 1/2, and there is a disjoint set of coordinates which does the same for Bob. It is easy to modify the example in a way that that it is not impossible for both of them to win.
		
		To see that in this case we can restrict  $L$ to be $K^c$ in \eqref{def:stbox}, note that for any sequence of tosses ${\bf x} \in A \stbbox B$ described above, we can take $K$  and $L$ to be the sets  of coordinates where heads and tails appear respectively, and then
		${{\bf x}} \in A_{1/2,K} \cap B_{1/2,L}$.
		$\qex$
	\end{example}

	\begin{example} [A three sided coin]  \label{ex:threeside}
		The goal of this example is to show that changing the definition of $A \stbbox B$ from allowing any disjoint $K$ and $L$, to always setting $L=K^c$, may result in a different set, unlike the case of $A \Box B$.
		Consider a three sided fair coin that shows Heads, Tails, or Sides. The coin is tossed 5 times, and similar to   Example \ref{ex:coin}, $A$ is the event that the first three tosses are heads, and $B$ that the last three are tails. Let {$(s,t)=(1/3,1/3)$}. Then if we take $K=\{1,2\}$ and $L=\{4,5\}$ we see that 
		${\bf x}= (HHSTT) \in A \stbbox B$.  However, if we only allowed $L=K^c$
		in the definition of $A \stbbox B$, then ${\bf x}$ would not be in this set.$\qex$
	\end{example}
	
	\begin{example}[Continuous percolation with positive-probability annihilation] \label{ex:contq} In Example \ref{ex:percann} we considered a graph model having randomly placed vertices where certain edges were annihilated with probability zero. In particular, if the points $X_i$ and $X_j$,  sampled independently from a continuous distribution, coincide, that is, if they are within distance $q=0$, then all their incident edges are annihilated.

		Generalizing this model, we annihilate all edges incident on any two distinct points that lie within distance $q \ge 0$ of each other. Assume that $X_1,\ldots,X_n$ are i.i.d. with a strictly positive density, bounded above by $c$ on $[0,1]^2$. Edges are formed as in Example \ref{ex:perc}, with annihilation as just described. When $q>0$, annihilation occurs with positive probability, and hence it is easy to see that when $A$ and $B$ are the crossing events given in Example \ref{ex:perc}, we have $A \eebbox B = \emptyset$. Indeed, similar to the reasoning in the proof of Proposition \ref{prop:1}, for each ${\bf x} \in A \eebbox B$ at least one of $|K|$ or $|L|$ cannot be greater than $n/2$. Taking this set to be $K$ without loss of generality, its complementary set contains enough coordinates so that the points it indexes can annihilate, with positive probability, all edges incident on points indexed by $K$.  
		$\qex$
	\end{example}
	
	The next proposition shows that in Example \ref{ex:contq}, if the thresholds $s$ and $t$ are sufficiently small relative to $q$, then $A \stbbox B$ is non-empty. To simplify matters, we now restrict to metrics $d$ in $\mathbb{R}^2$ that are induced by norms. Given a norm, let $\tau$ denote the area of the unit ball $\{x:d(0,x) \le 1\}$ in the plane. Then the area of balls of radius $q \ge 0$ in the metric space $([0,1]^2,d)$ is bounded by $\tau q^2$, with equality if the entire disk is in $[0,1]^2$.

	We use the following notation: for $K \subset [n]$ let  $\mathscr{A}_K$ denote the event that there is crossing from left to right corresponding to the event $A$, however now only the points $\{X_i: i \in K\}$ are taken into consideration, that is, edges are only formed between pairs of points indexed by $K$ having distance $\le 2r$, and only points indexed by $K$ may annihilate any of these edges. 
	The event	$A$ that there is crossing from left to right when all $n$ points are taken into consideration, as in Example \ref{ex:perc}, through edges that are not annihilated, coincides with $\mathscr{A}_{[n]}$. It is easy to see that for $n$ and $|K|$ large enough, $\mathscr{A}_K$ has positive probability. Similarly, for $L \subset [n]$ let $\mathscr{B}_L$ denote the corresponding event for crossing from bottom to top.
	\begin{proposition}\label{prop:emptybst} Let $A$ be the event that there is crossing from left to right as in Example \ref{ex:perc} through edges that are not annihilated, and let $B$ denote such a crossing from bottom to top. 
		Let
		${\bf x} = (x_1,\ldots,x_n) \in ([0,1]^2)^n$  be such that ${\bf x} \in \mathscr{A}_K$ for some $K \subset [n]$, and ${\bf x} \in \mathscr{B}_L$
		for a disjoint set $L \subset [n]$. 
		Then ${\bf x} \in A \stbbox B$ provided 
		$s \,\,\mbox{and}\,\, t$ 
		lie in the interval
		$[0, 1-\frac{c}{4}n^2 \tau q^2]$, and $q$ is small enough so that $\frac{c}{4}n^2 \tau q^2 <1$, where $c$ is an upper bound on the common density function of $X_1,\ldots,X_n$.
	\end{proposition}
	\proof   Consider $A$, reasoning for $B$ being similar, and let ${\bf x}$ and $K$ satisfy the conditions of the proposition. The number of points that make up the path corresponding to $A$ is at most $|K|$, and hence there is at most area $|K| \tau q^2$ into which a point having index in ${K^c}$ can land to cause annihilation, which happens with probability at most $c|K|\tau q^2$. Hence, letting $E_j$ be the event that the point indexed by $j \in {K^c}$ causes the annihilation of the edges incident to some point on the path, the conditional probability of the event $F$ that some point on the path is annihilated, given $\{{X}_i,\, i \in K\}$ satisfies
	\begin{multline} \label{eq:to.not.hang}
		P({\bf X \in}F\mid {\bf X}_K={\bf x}_K)=P\Big(\bigcup_{j \in {K^c}}E_j\mid {\bf X}_K={\bf x}_K\Big)  \\
		\le \sum_{j \in {K^c}} P(E_j\mid {\bf x}_K) \le (n-|K|)|K| c \tau q^2 \le 
		\frac{c}{4}n^2 \tau q^2 \le 1-s.
	\end{multline}
	Given ${\bf x} \in \mathscr{A}_K$, that is, a realization in which the points indexed by $K$ `already' form a path,  $A$ can only fail to occur if some edge incident to a point on this path is annihilated. Hence 
	$P({\bf X} \in {A}\mid {\bf {X}}_K={\bf x}_K) {\ge} P({\bf X \in}F^c\mid{\bf X}_K={\bf x}_K) \ge s$. Repeating the argument for $B$, the result follows in view of \eqref{def:firstst}.
	\bbox
	
	Another way to think about the situation considered in Proposition \ref{prop:emptybst} is to consider the case where points indexed by $K$ have been revealed, and, perhaps despite annihilation caused by the revealed points, form a path corresponding to $A$. We then ask about the conditional probability that the remaining, unobserved points will destroy this path by annihilation. It is easy to see that for  
	$q$ satisfying the condition of the propostion the probability that some unobserved point
	annihilates the edges of an observed one is at most $1-(1-|K|c \tau q^2)^{n-|K|} \le (n-|K|)|K| c \tau q^2$, agreeing with the bound in \eqref{eq:to.not.hang}.

	Note that $r$ of Example \ref{ex:perc} can be chosen so small that $\lceil{n/2}\rceil$ points may not suffice to form a path as required by the events $A$ and $B$. In this case disjoint $K$ and $L$, both of which must demonstrate paths, cannot exist. Nevertheless, in such a case $A \stbbox B$ may still be non-empty, as the coordinates of ${\bf x}$ revealed by $K$ and $L$ when considered only marginally may be sufficiently promising for the chances that the
	yet unrevealed coordinates will yield $A$, and the same for $B$.

	\section{Earlier Results on Continuous Spaces} \label{sec:GuRa}
	We conclude with a brief mention of an earlier BKR result for the
	special case of continuous spaces given in \cite{GuRa}; see also references therein. The Poisson Boolean percolation model consists of an inhomogeneous Poisson point process in some bounded region  of ${\mathbb{R}}^d$, each point having a ball of a random radius around it, independent of the point process. For events $A$ and $B$ defined by this process, $A \Box B$ comprises those configurations of process points and their associated radii for which the occurrence of $A$ and $B$ can respectively be claimed upon revealing all process points, and their associated balls, lying in two disjoint spatial regions.  In this setting, inequality \eqref{eq:box.inequality} is proved in a way that relies heavily on properties of the Poisson process. Consequently, it appears that this Poisson process case does not imply our results for a fixed number of points. 
	
	Conversely, even leaving the spatial structure in \cite{GuRa} aside, it seems that our results do not lead to results that involve a random number of points.  
	We can condition on the number of Poisson points and obtain $P(A \Box B\mid N) \le P(A\mid N)P(B\mid N)$ and take expectation with respect to $N$. However, the inequality in $E\{P(A\mid N)P(B\mid N)\} \le E\{P(A\mid N)\}E\{P(B\mid N)\} =P(A)P(B)$ holds only under very special conditions.

	\section{ Appendix on Measurability}\label{sec:meas}
	In this section we prove Propositions \ref{prop:partof4}, \ref{prop:LebegT} and \ref{prop:measure.st.sets}, the first two of which pertain to Polish spaces. \smallskip

	\noindent {\em Proof of Propositions \ref{prop:partof4} and \ref{prop:LebegT}:}
	Both proposition follow once we show  that if $A \in \bbs$ then $[A]_K$ is measurable with respect to the Borel sigma algebra $\bbs$ when completed relative to the product probability measure $P$ of Borel measures on a product of Polish spaces. For Polish spaces restricted to be subspaces of $\mathbb{R}$ this was shown 
	in Section 8 and Lemma 1 of \cite{AGH}.

	First, as in \cite{AGH}, with ${\rm Proj}_K$ denoting the projection of $x \in S$ on the coordinates in $K$, we may express the cylinders $[A]_K$ defined in \eqref{def:AGH.cyl.box1} as
	\bea \label{def:cylinder}
	[A]_K=\{ y \in S : y_K \in S_K \setminus {\rm Proj}_K(A^c) \}.
	\ena
	
	Letting $A$ be a set in the product Borel sigma algebra, the complement $A^c$ is Borel. By \cite{Cohn}, Proposition 8.4.4 (see also the text preceding) the projection $C={\rm Proj}_K(A^c)$ is analytic. 
	In particular, $C$ is universally measurable, and hence measurable in the completion of $P_K$. 
	The set $E=S_K \setminus {\rm Proj}_K(A^c)$, the complement of $C$ in $S_K$, is therefore also measurable in the completion of $P_K$.
	Letting $K=\{1,\ldots,k\}$ without loss of generality, note that we may explicitly write the ``pullback" $[A]_K$ of $E$ as
	\beas
	[A]_K=E \times \prod_{j=k+1}^n S_j.
	\enas
	As $E$ is measurable with respect to the completion of $P_K$, there exists $F$ and $G$ in ${\cal S}_K$ such that $F \subset E \subset G$ and $P_K(F)=P_K(G)$, yielding the inclusion
	\beas
	F \times \prod_{j=k+1}^n S_j \subset [A]_K \subset G\times \prod_{j=k+1}^n S_j.
	\enas
	As the sets on the left and right are elements of the product Borel sigma algebra $\bbs$, and have the same $P$ measure, $[A]_K$ is measurable with respect to the completion of $P$. This argument, as in \cite{AGH}, holds also for countable products.
	
	\noindent {\em Proof of Proposition \ref{prop:measure.st.sets}.}
	A proof that for ${\bf y} \in S_K$ and $E \in \mathbb S$ the set $\{{\bf v} \in S_{K^c}: {\langle {\bf y}, {\bf v}  \rangle_K} \in E\}$ is measurable is given in \cite{Fo} Proposition 2.34a.
	Next we show that  $P_{K^c}\left(\{{\bf v} \in S_{K^c}: {\langle {\bf y}, {\bf v}  \rangle_K} \in E\}\right)$ is measurable in ${\bf y} \in S_K$, and that it is a version of the conditional probability $P({\bf X} \in E|{\bf X}_K={\bf {y}})$.

	Fix some $K \subseteq [n]$ and $E \in \mathbb S$. For ${\bf y} \in S_K$ let $Q_{K^c}({\bf y};E):=P_{K^c}\{{\bf v} \in S_{K^c}: \langle {\bf y}, {\bf v}  \rangle_K \in E\}$. Dropping the subscript ${K^c}$ when not needed for clarity, that $Q({\bf y};E)$ is a measurable function of ${\bf y}$ follows directly from Tonelli's Theorem, see \cite{Fo}, 2.37b, as
	\bea \label{eq:Q.is.int}
	Q({\bf y};E) = \int_{S_{K^c}} \phi_E({\bf y}, {\bf v})P_{K^c}(d{\bf v}) \qmq{{upon setting} ${\phi}_E({\bf y},{\bf v})={\bf 1}(\langle {\bf y}, {\bf v} \rangle_K \in E)$.}
	\ena
	It is immediate that $Q({\bf y};\cdot)$ is a probability measure on $(S,{\mathbb S})$ for each ${\bf y} \in S_K$.
	
	We next show that $Q({\bf y};E)$ is a version of $P({\bf X} \in C \mid {\bf X}_K={\bf {y}})$, that is, following 4.18 in \cite{Brei}, that
	\bea \label{eq:DC}
	\int_{{D}}
	{Q({\bf y};E) P_K({d\bf y})}
	= P(\{{\bf X} \in E\} \cap \{{\bf X}_K \in D\})  {\qmq{for all $D \in \bbs_K$.}}
	\ena
	Using \eqref{eq:Q.is.int}, starting from the left hand side of \eqref{eq:DC}, we have, using $P=P_{K^c}P_K$ for the second equality, followed by Tonelli's Theorem,
	\begin{multline*}
		\int_D Q({\bf y};E)P_K(d{\bf y}) = \int_D \int_{S_{K^c}} \phi_E({\bf y}, {\bf v})P_{K^c}(d{\bf v})P_K(d{\bf y})
		\\
		={\int_D \int_{S_{K^c}} {\bf 1}(\langle {\bf y}, {\bf v} \rangle_K \in E)
			P_{K^c}(d{\bf v})P_K(d{\bf y})
			= \int {\bf 1}\left({\bf x} \in E, [{\bf x}]_K \in D \right)
			P(d{\bf x})}
		=P(\{{\bf X} \in E\} \cap \{{\bf X}_K \in D\}),
	\end{multline*}
	which is \eqref{eq:DC}. 
	
	Next we show measurability of $E_{r,{K}}$ in \eqref{def:uv.cond}, $E_r$ of \eqref{eq:nowinflate} and $A \stbbox B$ of \eqref{def:stbox}; as by \eqref{def:uv.cond} the sets $E_{r,K}$ for $K=[n]$ are clearly measurable, we restrict to the case $|K| \le n-1$.

	Letting again ${\rm Proj}_K:S \rightarrow S_K$ be the projection that sends ${\bf x}$ to ${\bf x}_K$, returning to \eqref{def:uv.cond} and considering ${\bf y} \mapsto Q_{K^c}({\bf y};E)$  as a function only of ${\bf y} \in S_K$ for fixed $K$ and $E$, we see that $E_{r,K}={\rm Proj}_K^{-1}(Q_{K^c}^{-1}([r,1];{E}))$, the inverse image of a composition of measurable mappings of a measurable set, and is therefore measurable.
	By  \eqref{eq:nowinflate} and  \eqref{def:stbox} it readily follows that $E_r$ and $A \stbbox B$ are measurable.

	Lastly, we show that the replacement in \eqref{def:uv.cond} of our 
	version of the conditional probability by another changes $E_{r,K}$ by a $P$-null set. Again, as $E_{r,[n]}$ given in \eqref{def:uv.cond} does not depend on the version, we restrict to the case  $|K| \le n-1$. Letting $G$ be the subset of $S_K$ where the two versions differ we have $P_K(G)=0$, and without loss of generality, let $K=\{1,\ldots,k\}$. With $E^i_{r,K}, i =1,2$ denoting the corresponding sets under these versions, the symmetric difference $E^1_{r,K} \triangle  E^2_{r,K}$ is contained in $G \times S_{K^c}$, a $P$-null set. Hence the symmetric difference has $P$ measure zero.  \qed

	\noindent {\bf Acknowledgements:} We are deeply indebted to an anonymous referee for a very careful reading of two versions of this paper. The penetrating comments provided enlightened us on various measurability issues and other important, subtle points.
	We thank Mathew Penrose for a useful discussion and for providing some relevant references. 
	
	The work of the first author was partially supported by NSA grant H98230-15-1-0250. The second author would like to thank the Isaac Newton Institute for Mathematical Sciences, Cambridge, for support and hospitality during the programme Data Linkage and Anonymisation where part of the work on this paper was undertaken,
	supported by EPSRC grant no EP/K032208/1.

\end{document}